\documentclass[leqno,11pt]{amsart}

\usepackage{latexsym,amssymb,amsthm,amsmath,amscd}

\theoremstyle{plain}
\newtheorem{theorem}{Theorem}

\numberwithin{equation}{section}

\theoremstyle{remark}
\newtheorem{remark}{Remark}[section]

 \numberwithin{equation}{section}

\def\<{\left < }
\def\>{\right >}
\def\({\left ( }
\def\){\right )}

\def\n2{\left[{n\over2}\right]}

\begin{document}

\title[Conformal mappings and first eigenvalue of Laplacian]{Conformal mappings and first eigenvalue of Laplacian on surfaces}
\author[B.-Y. Chen]{Bang-Yen Chen}
\address{Chen: Department of Mathematics\\
	Michigan State University \\East Lansing, Michigan 48824--1027\\ U.S.A.}
\email{bychen@math.msu.edu}

\begin{abstract}   In this note we give a simple relation between conformal mapping and the first eigenvalue of Laplacian for surfaces in Euclidean spaces.
\end{abstract}

\thanks{This article was published in Bulletin of the Institute of Mathematics, Academia Sinica, {\bf 7} (1979), 395--400.} 
\maketitle

\section{Statement of Main Theorem.} 

Let $M$ be a compact Riemannian surface and $\Delta$ the Laplace-Beltrami operator acting on differential functions $C^\infty(M)$ on $M$. It is known that $\Delta$ is an elliptic operator. The operator $\Delta$ has an infinite sequence
\begin{align} 0=\lambda_0<\lambda_1<\lambda_2<\cdots<\lambda_p <\cdots \nearrow \infty \end{align}
of eigenvalues. 

Let $V_p=\{f\in C^\infty(M):\Delta f=\lambda_[  f\}$ be the eigenspace with eigenvalue $\lambda_p$. Then the dimension of $V_p$ is finite, which is called the multiplicity of $\lambda_p$. 
Let $x:M\to E^m$ be an immersion of a surface $M$ in $E^m$. Then the Euclidean metric of $E^m$ induces a Riemannian metric on $M$. In this paper we shall consider only the induced metric on $M$. As in [4], we shall call an immersion $x:M\to E^m$  of order $p$ if all coordinate functions of $x=(x_1,\ldots,x_m)$ are in $V_p$, where $x_1,\ldots,x_m$ are the Euclidean coordinates of $x$. 

In the following we shall denote by $\lambda_1(x)$ and $A(x)$ respectively the first eigenvalue $\lambda_1$ and the area of $M$ with respect to the immersion $x$ when it is necessary.

In this paper, we shall prove the following {\it conformal inequality} for $\lambda_1$.

\begin{theorem} Let $x:M\to E^m$ be an imbedding of order 1 from a compact surface $M$ in $E^m$. If $\varphi$ is a conformal mapping of $E^m$ with $A(x)=A(\varphi\circ x)$, then we have
\begin{align} \lambda_1(\varphi\circ x)\leq \lambda_1(x).\end{align}
The equality holds when and only when $\varphi$ is a rigid motion of $E^m$.
\end{theorem}

If $c$ is a nonzero constant, then we have $$\lambda_1(cx)=\frac{\lambda_1(x)}{c^2},\;\;  A(cx)=c^2A(x)$$ for the similarity transformation $cx$. Thus the assumption on the area in Theorem 1 is necessary. Moreover the assumption is generic in the sense that if $A(x)\ne A(\varphi\circ x)$, then by choosing a suitable similarity transformation $\psi$ of $E^m$ we have $A(x)=A(\psi\circ \varphi\circ x)$.

It seems to the author that inequality (1.2) is the only conformal inequality we know so far for spectra. Some applications of Theorem 1 will be given in the last section. A typical example reads as follows: {\it If $M$ is a cyclide of Dupin given by an inversion of an anchor ring in $E^3$ with circles of radii $a$ and $b$ satisfy $a/b=\sqrt{2}$, then $\lambda_1<4\pi^2/A$.}

\section{Proof of Theorem 1.}

Let $x:M\to E^m$ be an imbedding of a compact surface $M$ in $E^m$. Without loss of generality we may choose the center of gravity as the origin of $E^m$. Let $x=(x_1,\ldots,x_m)$ be the Euclidean coordinates of $E^m$. Then we have $\int_M x_idV=0$. The minimal principle [1] then implies
\begin{align} \int_M |dx_i|^2dV\geq  \lambda_1\int_M (x_i)^2dV,\;\; i=1,\ldots,m.\end{align}
The equality holds if and only if each $x_i$ is in $V_1$. On the other hand, since $|dx|^2=\sum_{i=1}^m |dx_i|^2 =2$, (2.1) gives
\begin{align}2A(x)\geq \lambda_1(x)\int_M |x|^2dV.\end{align}
Let $H=\frac{1}{2}\, {\rm trace}\,\sigma$ be the mean curvature vector of $x$, $\sigma$ the second fundamental form of $x$. Then we have [6]
\begin{align}  A(x)+\int_M (x\cdot H)dV=0.\end{align}
Thus (2.2), (2.3) and the Schwartz inequality imply
\begin{align}\notag  2A\int_M |H|^2dV&\geq \lambda_1\(\int_M (|x|\,| H|)dV\)^2
\\&\notag \lambda_1 \(\int_M (x\cdot H)dV\)^2\geq \lambda_1 A^2.\end{align}
Consequently, we obtain the following result of Reilly [7]
\begin{align} \int_M |H|^2dV\geq  \frac{\lambda_1(x)}{2}A(x).\end{align}
The equality holds if and only if $x-a$ is of order 1 for some vector $a$ in $E^m$.

In the following, we denote by $TMC(x)$ the total mean curvature of $x$, i.e.,
$$TMC(x)=\int_M  |H|^2dV.$$
Suppose that $x:M\to E^m$ is an imbedding of order 1. Then we have
\begin{align} TMC(x)=\frac{\lambda_1(x)}{2}A(x).\end{align}

If $\varphi$ is a conformal mapping of $E^m$ into $E^m$, then we have [3]
\begin{align} TMC(\varphi\circ x)=TMC(x).\end{align}
Combining this with (2.4) and (2.5), we find
\begin{align} \lambda_1(\varphi\circ x)A(\varphi\circ x)\leq \lambda_1(x)A(x).\end{align}
The equality holds if and only if $\varphi\circ x-b$ is of order 1 for some vector $b$ in $E^m$. In particular, if $A(x)=A(\varphi\circ x)$, (2.7) gives
\begin{align} \lambda_1(\varphi\circ x)\leq\lambda_1(x).\end{align}

If the equality of (2.8) holds, the $\varphi\circ x-b$ is also of order 1 for some vector $b$ in $E^m$. By using a translation on $E^m$, we may also assume that the center of gravity of $\varphi\circ x$ is the origin too. In this case, $b=0$ and $\varphi\circ x$ is of order 1. Consequently, we have
\begin{align} \Delta x=\lambda_1x,\;\;\bar\Delta (\varphi\circ x)=\lambda_1(\varphi\circ x),\end{align}
where $\lambda_1=\lambda_1(x)=\lambda_1(\varphi\circ x)$ and $\Delta$ and $\tilde\Delta$ are the Laplace-Beltrami operators on $M$ with respect to $x$ and $\varphi\circ x$, respectively. 
From (2.9) and a theorem of Takahashi [8], we see that $M$ is imbedded by $x$ and $\varphi\circ x$ into the same hypersphere $S^{m-1}(r)$ of radius $r=\sqrt{2/\lambda_1}$
 as minimal surfaces.

Now, by a result of Haantjes [5], we know that a conformal mapping on $E^m$ is generated by translations, rotations, homothetic transformations and inversions centered at a fixed point. Since the centers of gravity of $x$ and $\varphi\circ x$ are assumed to be at the same point 0, the conformal mapping $\varphi$ is free of translation. Moreover, since $M$ is imbedded both by $x$ and $\varphi\circ x$ into the same hypersphere $S^{m-1}(r)$, $\varphi$ is free of homothetic transformations (except the identity transformation). On the other hand, inversions centered at 0 are given in the following form:
$$\bar x=\frac{c^2}{(x\cdot x)}x$$
for nonzero constants $c$. Thus inversions centered at 0 always carry a hypersphere of radius $r$ into a hypersphere of radius $c^4/r^2$. In our case, since both surfaces given by $x$ and $\varphi\circ x$ lie in the same hypersphere $S^{m-1}(r)$, $\varphi$ is free of inversions too. Consequently, $\varphi$ is given only by a rotation. Conversely, because the area and the spectrum of a surface are invariant under rigid motion (generated by translations and rotations), if $\varphi$ is a rigid motion, the equality of (1.2) holds.

\begin{remark} Theorem 1 shows that the estimates on total mean curvature for surfaces in $E^m$ given in [2,7] are weak in general.
\end{remark}

\section{Applications.}

In this section we shall give the following applications of Theorem 1.

Let $S^1(1)$ be the unit circle in a plane $E^2$. Then the product surface $T^2=S^1(1)\times S^1(1)$ is a flat surface in $E^4$ with area $A=4\pi^2$ and $\lambda_1=1$ [1]. A surface in $E^m\supset E^4$ is called a {\it conformal Clifford torus} if it is the image of the Clifford torus under a conformal mapping of $E^m$. The anchor ring in $E^3$ given by
$$ \Big((\sqrt{2}+\cos u)a \cos v, (\sqrt{2}+\cos u)a\sin v, a\sin v\Big)$$
$$0\leq u<2\pi,\;\; 0\leq v<2\pi,$$
is among the class of conformal Clifford tori. It is easy to see that the Clifford torus in $E^4$ is order 1. There exists no conformal Clifford torus of order 1 in $E^3$. 

Theorem 1 implies the following.

\begin{theorem} Let $M$ be a conformal Clifford torus in $E^m$ $(m\geq 3)$ with area $4\pi^2$, Then we have
\begin{align} \lambda_1\leq 1.\end{align}
The equality holds if and only if $M$ is a Clifford torus.
\end{theorem}

Let $(x,y,z)$ be the Euclidean coordinates of $E^3$ and $(u^1,\ldots,u^5)$ be the Euclidean coordinates of $E^5$. We consider the mapping defined by
$$u^1=\frac{1}{3} yz, \;\; u^2=\frac{1}{3}zx,\;\; u^3=\frac{1}{3}xy,$$
$$u^4=\frac{1}{6}(x^2-y^2),\;\; u^5=\frac{1}{6\sqrt{3}}(x^2+y^2-2z^2).$$
This defines an isometric immersions of $S^2(1)$ into $S^5(1/\sqrt{3})$ as a minimal surface. Two points $(x,y,z)$ and $(-x,-y,-z)$ of $S^2(1)$ are mapped into the same point of $S^4(1/\sqrt{3})$ and this mapping defines an imbedding of the real projective plane in  $S^4(1/\sqrt{3})\subset E^5$. This real projective plane imbedded in $E^5$ is called the {\it Veronese surface}. It is known that Veronese surface satisfies $A=2\pi$ and $\lambda_1=6$ [1]. A surface in $E^m$ $(m\geq 5)$ is called a {\it conformal Veronese surface} if it is the image of the Veronese surfaces under a conformal mapping of $E^m$. There is no conformal Veronese surface of order 1 in $E^4$. 

From Theorem 1 we have

\begin{theorem} Let $M$ be a conformal Veronese surface with area $2\pi$. Then we have
\begin{align} \lambda_1\leq 6.\end{align}
The equality holds if and only if $M$ is a Veronese surface.
\end{theorem}

\end{document}